\documentclass{article}
\usepackage[utf8]{inputenc}
\usepackage{geometry}
\usepackage{amsmath}
\usepackage{amssymb}
\usepackage{amsthm}
\usepackage{tgbonum}
\usepackage{cite}
\usepackage{tikz}
\usetikzlibrary{
  knots,
  hobby,
  decorations.pathreplacing,
  decorations.markings,
  shapes.geometric
}

\theoremstyle{plain}
\newtheorem{theorem}{Theorem}
\newtheorem{lemma}{Lemma}
\theoremstyle{definition} 
\newtheorem{definition}{Definition}

\newcommand{\make}[1]{\; 
	\begin{tikzpicture}
		[baseline=-\dimexpr\fontdimen22\textfont2\relax] #1
	\end{tikzpicture} \; 
}
\newcommand{\makeO}{ \make{ 
	\filldraw[color=gray, fill=none, thick] circle (0.25);
}}
\newcommand{\makeXN}{\make{ 
    \draw[color=gray,thick] (-0.3,0.3) -- (0.3,-0.3);
    \draw[color=gray,thick] (-0.3,-0.3) -- (-0.1,-0.1);
    \draw[color=gray,thick] (0.1,0.1) -- (0.3,0.3);
  }}
\newcommand{\makeXND}{\make{ 
    \draw[->,color=gray,thick] (-0.3,0.3) -- (0.3,-0.3);
    \draw[<-, color=gray,thick] (-0.3,-0.3) -- (-0.1,-0.1);
    \draw[color=gray,thick] (0.1,0.1) -- (0.3,0.3);
  }}

\newcommand{\makeXPD}{\make{ 
   \draw[color=gray,thick] (-0.3, 0.3) -- (-0.1, 0.1);
   \draw[->, color=gray,thick] (0.1,-0.1) -- (0.3,-0.3);
   \draw[<-,color=gray,thick] (-0.3,-0.3) -- (0.3,0.3);
  }}
\newcommand{\makeXV}{\make{ 
    \draw[color=gray,thick] (-0.3,0.3) -- (0.3,-0.3);
    \draw[color=gray,thick] (-0.3,-0.3) -- (0.3,0.3);
  }}
\newcommand{\makeC}{\make{ 
    \draw[color=gray,thick] (-0.3,0.3) .. controls (0,-0.025) .. (0.3,0.3);
    \draw[color=gray,thick] (-0.3,-0.3) .. controls (0,0.025) .. (0.3,-0.3);
  }}
\newcommand{\makeTT}{\make{
    \draw[color=gray,thick] (-0.3,-0.3) .. controls (0.025,0) .. (-0.3,0.3);
    \draw[color=gray,thick] (0.3,-0.3) .. controls (-0.025,0) .. (0.3,0.3);
  }}
\newcommand{\makePun}{\make{
	\draw[color=gray,thick] (-0.3,0.3) -- (-0.3,-0.3)--(0.3,-0.3)--(0.3,0.3)--(-0.3,0.3);
	\node at (0,0) {$L$};
	\draw[color=gray,thick] (0.3, 0.2) -- (0.5, 0.07);
	\draw[color=gray,thick] (0.7,-0.07) -- (0.9,-0.2);
	\draw[color=gray,thick,decoration={markings, mark=at position 0.2 with {\arrow{<}}},postaction={decorate}] (0.3,-0.2) -- (0.9,0.2) arc(90:-90:.2)--(.9,-.2) ;
  }}
\newcommand{\makeFunn}{\make{
	\draw[color=gray,thick] (-0.3,0.3) -- (-0.3,-0.3)--(0.3,-0.3)--(0.3,0.3)--(-0.3,0.3);
	\node at (0,0) {$L$};
	\draw[color=gray,thick, decoration={markings, mark=at position 0.5 with {\arrow{>}}}, postaction={decorate}] (0.3,0.2) arc(90:-90:.2) ;
  }}
\newcommand{\makeVun}{\make{
	\draw[color=gray,thick] (-0.3,0.3) -- (-0.3,-0.3)--(0.3,-0.3)--(0.3,0.3)--(-0.3,0.3);
	\node at (0,0) {$L$};
	\draw[color=gray,thick] (0.3, 0.2) -- (0.9, -0.2);
	\draw[color=gray,thick] (0.3,-0.2) -- (0.9,0.2) arc(90:-90:.2)--(.9,-.2) ;
  }}
\newcommand{\makeNun}{\make{
	\draw[color=gray,thick] (-0.3,0.3) -- (-0.3,-0.3)--(0.3,-0.3)--(0.3,0.3)--(-0.3,0.3);
	\node at (0,0) {$L$};
	\draw[color=gray,thick,decoration={markings, mark=at position 0.425 with {\arrow{>}}}, postaction={decorate}] (0.3, 0.2) -- (0.9, -0.2) arc(-90:90:.2)--(.9,.2) ;
	\draw[color=gray,thick] (0.3,-0.2) -- (0.5,-0.07); 
     \draw[color=gray,thick] (0.7,0.07) -- (0.9,0.2); 
  }}
\newcommand{\makeFun}{\make{
	\draw[color=gray,thick] (-0.3,0.3) -- (-0.3,-0.3)--(0.3,-0.3)--(0.3,0.3)--(-0.3,0.3);
	\node at (0,0) {$L$};
	\draw[color=gray,thick] (0.3,0.2) arc(90:-90:.2) ;
  }}
\newcommand{\makeWun}{\make{
	\draw[color=gray,thick] (-0.3,0.3) -- (-0.3,-0.3)--(0.3,-0.3)--(0.3,0.3)--(-0.3,0.3);
	\node at (0,0) {$L$};
	\draw[color=gray,thick] (0.3,0.2) arc(90:-90:.2) ;
	\draw[color=gray,thick] (0.45,.05) -- (0.55,-.05);   
  }}
  
\newcommand{\makeRRhalf}{\make{
	\draw[color=gray,thick, ->] (0.2,0.4) .. controls (-0.25,0) .. (0.2,-0.4); 
		\draw [fill, color=white] (0.01,0.21) circle (.052); 
		\draw [fill, color=white] (0.01,-0.21) circle (.052);
	\draw[color=gray,thick, ->] (-0.2,0.4) .. controls (0.25,0) .. (-0.2,-0.4);	
  }}

\newcommand{\makeRp}{\make{
	\draw[thick, ->] (-.9,0.4) arc(40:-40:0.625);
	\node at (-.35,0) {$\leftrightarrow$}; \node [above] at (-.35,0) {\tiny{(R1a)}};
	\draw[color=black,thick] (0, 0.4) to [out=-65, in=120] (0.2,0) to [out=-60, in=180] (0.5,-0.2) to [out=0, in=270](.75,0);
			\draw [fill, color=white] (0.2,0) circle (.052); 
	\draw[color=black,thick, <-] (0, -0.4) to [out=65, in=240] (0.2,0) to [out=60, in=180] (0.5,0.2) to [out=0, in=90] (.75,0);
	
	\draw [fill, color=white] (0,0.6) circle (.052); 
  }}
\newcommand{\makeRv}{\make{
	\draw[thick, ->] (-.8,0.4) arc(40:-40:0.625) ;
	\node at (-.3,0) {$\leftrightarrow$}; \node [above] at (-.3,0) {\tiny{(V1)}};
	\draw[color=black,thick] (0, 0.4) to [out=-65, in=120] (0.2,0) to [out=-60, in=180] (0.5,-0.2) to [out=0, in=270](.75,0);
	\draw[color=black,thick, <-] (0, -0.4) to [out=65, in=240] (0.2,0) to [out=60, in=180] (0.5,0.2) to [out=0, in=90] (.75,0);
	\draw [fill, color=white] (0,0.6) circle (.052); 
  }}
\newcommand{\makeRn}{\make{
    \draw[thick, ->] (-.9,0.4) arc(40:-40:0.625);
	\node at (-.35,0) {$\leftrightarrow$}; \node [above] at (-.35,0) {\tiny{(R1b)}};
	\draw[color=black,thick, <-] (0, -0.4) to [out=65, in=240] (0.2,0) to [out=60, in=180] (0.5,0.2) to [out=0, in=90] (.75,0);
			\draw [fill, color=white] (0.2,0) circle (.052); 
	\draw[color=black,thick] (0, 0.4) to [out=-65, in=120] (0.2,0) to [out=-60, in=180] (0.5,-0.2) to [out=0, in=270](.75,0);
  }}  
\newcommand{\makeVV}{\make{
	\draw[thick, ->] (-0.2,0.4) .. controls (0.25,0) .. (-0.2,-0.4);
    \draw[thick, ->] (0.2,0.4) .. controls (-0.25,0) .. (0.2,-0.4);
	
	\node at (0.5,0) {$\leftrightarrow$}; \node [above] at (.5,0) {\tiny{(V2)}};
	\draw[thick, ->] (.8,0.4) arc(40:-40:0.625) ;
	\draw[thick, ->] (1.25,0.4) arc(140:220:.625) ; 
  }}
\newcommand{\makeRR}{\make{
	\draw[color=black,thick, ->] (0.2,0.4) .. controls (-0.25,0) .. (0.2,-0.4); 
		\draw [fill, color=white] (0.01,0.21) circle (.052); 
		\draw [fill, color=white] (0.01,-0.21) circle (.052);
	\draw[color=black,thick, ->] (-0.2,0.4) .. controls (0.25,0) .. (-0.2,-0.4);

	\node at (0.5,0) {$\leftrightarrow$}; \node [above] at (.5,0) {\tiny{(R2)}};
	\draw[color=black,thick, ->] (.8,0.4) arc(40:-40:0.625) ;
	\draw[color=black,thick, ->] (1.25,0.4) arc(140:220:.625) ; 
	
  }}

\newcommand{\makeVVV}{\make {
	\draw[color=black, thick, ->] (-.4, .4) to [out=-80, in=150]  (-0.2,0.2) to [out=-30, in=150]  (0.15,0)  to [out=-30, in=80] (.4, -.4); 
	\draw[color=black, thick, ->] (0, .4)  to [out=225, in=60]  (-0.2,0.2)  to [out=250, in=110]  (-0.2,-0.2)  to [out=-60, in=135]  (0, -.4); 
	\draw[color=black, thick, ->] (0.4, 0.4) to [out=-80, in=30]  (0.15,0) to [out=210, in=30]  (-0.2,-0.2)  to [out=210, in=80] (-.4, -.4); 
	
	\draw[color=black, thick, ->] (1.1, .4) to [out=-80, in=150]  (1.35,0) to [out=-30, in=150]  (1.7,-0.2) to [out=-30, in=100] (1.9, -.4); 
	\draw[color=black, thick, ->] (1.5, .4)  to [out=-45, in=120]  (1.7,0.2)  to [out=-70, in=70]  (1.7,-0.2)  to [out=-120, in=45]  (1.5, -.4); 
	\draw[color=black, thick, ->] (1.9, 0.4) to [out=-80, in=30]  (1.7,0.2)  to [out=210, in=30] (1.35,0)  to [out=210, in=100] (1.1, -.4);	
	\node at (.75,0) {$\leftrightarrow$}; \node [above] at (.75,0) {\tiny{(V3)}};
  }}  
\newcommand{\makeRRR}{\make{
	\draw[color=black, thick, ->] (-.4, .4) to [out=-80, in=150]  (-0.2,0.2) to [out=-30, in=150]  (0.2,0)  to [out=-30, in=80] (.4, -.4);
			\draw [fill, color=white] (-0.2,0.2) circle (.05);
			\draw [fill, color=white] (0.175,0.01) circle (.05);
	\draw[color=black, thick, ->] (0, .4)  to [out=225, in=60]  (-0.2,0.2)  to [out=250, in=110]  (-0.2,-0.2)  to [out=-60, in=135]  (0, -.4);
			\draw [fill, color=white] (-0.2,-0.2) circle (.05);
	\draw[color=black, thick,->] (0.4, 0.4) to [out=-80, in=30]  (0.15,0) to [out=210, in=30]  (-0.2,-0.2)  to [out=210, in=80] (-.4, -.4);
	
	\draw[color=black, thick, ->] (1.1, .4) to [out=-80, in=150]  (1.35,0) to [out=-30, in=150]  (1.7,-0.2) to [out=-30, in=100] (1.9, -.4);
    		\draw [fill, color=white] (1.7,-0.2) circle (.05);
	\draw[color=black, thick, ->] (1.5, .4)  to [out=-45, in=120]  (1.7,0.2)  to [out=-70, in=70]  (1.7,-0.2)  to [out=-120, in=45]  (1.5, -.4);
			\draw [fill, color=white] (1.7,0.2) circle (.05);
			\draw [fill, color=white] (1.35,0) circle (.05);
	\draw[color=black, thick, ->] (1.9, 0.4) to [out=-80, in=30]  (1.7,0.2)  to [out=210, in=30] (1.35,0)  to [out=210, in=100] (1.1, -.4);		
	\node at (.75,0) {$\leftrightarrow$}; \node [above] at (.75,0) {\tiny{(R3)}};
  }}    
  
\newcommand{\makeW}{\make{
	\draw[color=black, thick, ->] (-.4, .4) to [out=-80, in=150]  (-0.2,0.2) to [out=-30, in=150]  (0.2,0)  to [out=-30, in=80] (.4, -.4);
			\draw [fill, color=white] (0.175,0.01) circle (.05);
	\draw[color=black, thick, ->] (0, .4)  to [out=225, in=60]  (-0.2,0.2)  to [out=250, in=110]  (-0.2,-0.2)  to [out=-60, in=135]  (0, -.4);
			\draw [fill, color=white] (-0.2,-0.2) circle (.05);
	\draw[color=black, thick, ->] (0.4, 0.4) to [out=-80, in=30]  (0.15,0) to [out=210, in=30]  (-0.2,-0.2)  to [out=210, in=80] (-.4, -.4);
	
	\draw[color=black, thick, ->] (1.5, .4)  to [out=-45, in=120]  (1.7,0.2)  to [out=-70, in=70]  (1.7,-0.2)  to [out=-120, in=45]  (1.5, -.4);	
	\draw[color=black, thick, ->] (1.1, .4) to [out=-80, in=150]  (1.35,0) to [out=-30, in=150]  (1.7,-0.2) to [out=-30, in=100] (1.9, -.4);
			\draw [fill, color=white] (1.7,0.2) circle (.05);
			\draw [fill, color=white] (1.35,0) circle (.05);
	\draw[color=black, thick, ->] (1.9, 0.4) to [out=-80, in=30]  (1.7,0.2)  to [out=210, in=30] (1.35,0)  to [out=210, in=100] (1.1, -.4);	
	\node at (.75,0) {$\leftrightarrow$}; \node [above] at (.75,0) {\tiny{(F1)}};
		\draw [fill, color=white] (0,0.6) circle (.052); 
  }}

\newcommand{\makeM}{\make{

\draw[color=black, thick, ->] (-.4, .4) to [out=-80, in=150]  (-0.2,0.2) to [out=-30, in=150]  (0.2,0)  to [out=-30, in=80] (.4, -.4);

\draw [fill, color=white] (-0.2,0.2) circle (.05); 

\draw[color=black, thick, ->] (0, .4)  to [out=225, in=60]  (-0.2,0.2)  to [out=250, in=110]  (-0.2,-0.2)  to [out=-60, in=135]  (0, -.4); 

\draw[color=black, thick, ->] (0.4, 0.4) to [out=-80, in=30]  (0.15,0) to [out=210, in=30]  (-0.2,-0.2)  to [out=210, in=80] (-.4, -.4);

\draw[color=black, thick, ->] (1.1, .4) to [out=-80, in=150]  (1.35,0) to [out=-30, in=150]  (1.7,-0.2) to [out=-30, in=100] (1.9, -.4); 

\draw [fill, color=white] (1.7,-0.2) circle (.05);

\draw[color=black, thick, ->] (1.5, .4)  to [out=-45, in=120]  (1.7,0.2)  to [out=-70, in=70]  (1.7,-0.2)  to [out=-120, in=45]  (1.5, -.4);

\draw[color=black, thick, ->] (1.9, 0.4) to [out=-80, in=30]  (1.7,0.2)  to [out=210, in=30] (1.35,0)  to [out=210, in=100] (1.1, -.4); 

\node at (.75,0) {$\leftrightarrow$}; \node [above] at (.75,0) {\tiny{(M)}};

\draw [fill, color=white] (0,0.6) circle (.052); 

  }}  
    
\newcommand{\makeT}{\make{
	\draw[color=black,thick, ->] (-0.2,0.4) -- (-0.2,-0.4);
			\draw[color=black,thick] (-0.275,.05) -- (-0.125,-0.05); 
			\draw[color=black,thick] (-0.275,.16) -- (-0.125,0.06);
	\node at (0.2,0) {$\leftrightarrow$}; \node [above] at (.2,0) {\tiny{(T1)}};
	\draw[color=black,thick, ->] (.6,0.4) -- (.6,-0.4);

	\draw [fill, color=white] (0,0.6) circle (.052); 
		\draw [fill, color=white] (-1.50,0.05) circle (.052); 
  }}
\newcommand{\makeTt}{\make{
	\draw[color=black,thick, ->] (-0.2,0.4) to [out=-90, in=90] (-0.2,0.05) to [out=-60, in=135] (0,-0.15)to [out=-45, in=120]  (0.2,-0.4);
	\draw[color=black,thick, ->] (0.2,0.4)  to [out=-90, in=90] (0.2,0.05) to [out=-120, in=45] (0,-0.15)to [out=225, in=60]  (-0.2,-0.4);
			\draw[color=black,thick] (-0.275,.18) -- (-0.125,0.08);
			
	\node at (0.55,0) {$\leftrightarrow$}; \node [above] at (.55,0) {\tiny{(T2)}};
	\draw[color=black,thick, ->] (.90,0.4) to [out=-60, in=135] (1.1,0.15) to [out=-45, in=120] (1.3,-0.05)to [out=-90, in=90]  (1.3,-0.4);
	\draw[color=black,thick, ->] (1.3,0.4)  to [out=240, in=45] (1.1,0.15) to [out=225, in=60]  (0.9,-0.05) to [out=-90, in=90]  (0.9,-0.4);	
			\draw[color=black,thick] (1.225,-0.08) -- (1.375,-0.18);
  }}  
\newcommand{\makeTTT}{\make{
	\draw[color=black,thick, ->] (-0.2,0.4) to [out=-90, in=90] (-0.2,0.05) to [out=-60, in=135] (0,-0.15)to [out=-45, in=120]  (0.2,-0.4);
			\draw[color=black,thick] (-0.275,.18) -- (-0.125,0.08);
			\draw [fill, color=white] (0,-0.15) circle (.052); 	
	\draw[color=black,thick, ->] (0.2,0.4)  to [out=-90, in=90] (0.2,0.05) to [out=-120, in=45] (0,-0.15)to [out=225, in=60]  (-0.2,-0.4);
	
	\node at (0.55,0) {$\leftrightarrow$}; \node [above] at (.55,0) {\tiny{(T3)}}; 
	\draw[color=black,thick, ->] (.90,0.4) to [out=-60, in=135] (1.1,0.15) to [out=-45, in=120] (1.3,-0.05)to [out=-90, in=90]  (1.3,-0.4);
			\draw[color=black,thick] (1.225,-0.08) -- (1.375,-0.18);
			\draw [fill, color=white] (1.1,0.15) circle (.052); 
	\draw[color=black,thick, ->] (1.3,0.4)  to [out=240, in=45] (1.1,0.15) to [out=225, in=60]  (0.9,-0.05) to [out=-90, in=90]  (0.9,-0.4);	
  }} 
\newcommand{\makeTTTT}{\make{
	\draw[color=black,thick, ->] (-0.2,0.4) to [out=-90, in=90] (-0.2,0.05) to [out=-60, in=135] (0,-0.15)to [out=-45, in=120]  (0.2,-0.4);
			\draw[color=black,thick] (0.125,.18) -- (0.275,0.08);
			\draw [fill, color=white] (0,-0.15) circle (.052); 	
	\draw[color=black,thick, ->] (0.2,0.4)  to [out=-90, in=90] (0.2,0.05) to [out=-120, in=45] (0,-0.15)to [out=225, in=60]  (-0.2,-0.4);
	
	\node at (0.55,0) {$\leftrightarrow$}; \node [above] at (.55,0) {\tiny{(T4)}};
   	\draw[thick, color=black, ->] (0.9,0.6) .. controls (2.10,0.150) and (.10,0.05) .. (1.3,-0.4) .. controls (1.350,-0.450) .. (1.37,-0.6) ;
			\draw [fill, color=white] (1.1,0.1) circle (.052); 
  	\draw[thick, color=black, ->] (1.3,0.6) .. controls (.10,0.150) and (2.10,0.05) .. (0.9,-0.4).. controls (0.850,-0.450).. (0.83,-0.6) ;
   			\draw[color=black,thick] (0.9,-0.325) -- (1.05,-0.425);
   			
   	\draw [fill, color=white] (0,0.6) circle (.052); 
  }} 
  
\newcommand{\makecB}{\make{
	\draw[color=gray] circle (1);
	\draw[color=black,thick] (-0.7,0.7) to [out=130, in=190] (-.5,1.3) to [out=10, in=90](0,1);
	\draw[color=black,thick] (-0.7,-0.7) to [out=240, in=180] (-.13,-1.3) to [out=0, in=270] (0.7,-0.7);
	\draw[color=black,thick] (0.7,0.7) to [out=30, in=100] (1.35,0.2) to [out=280, in=70] (1.3,-0.75) to [out=245, in=5](.56,-1.3) to [out=190, in=-60](0,-1);
  }}
\newcommand{\makeCases}{\make{
	\draw[color=black, thick] circle (1);
	\node [left] at (-.7,.7) {1}; \draw [fill] (-.7,.7) circle (.05);
	\node [above] at (0,1) {2}; \draw [fill] (0,1) circle (.05);
	\node [right] at (.7,.7) {3}; \draw [fill] (.7,.7) circle (.05);
	\node [right] at (.7,-.7) {4}; \draw [fill] (.7,-.7) circle (.05);
	\node [below] at (0,-1) {5}; \draw [fill] (0,-1) circle (.05);
	\node [left] at (-.7,-.7) {6}; \draw [fill] (-.7,-.7) circle (.05);
  }}
\newcommand{\makeInA}{\make{ 
	\draw[color= gray] circle (1);] 
	\draw[color=black,thick] (-.7,.7) --(0,-1); 
	\draw[color=black,thick] (0, 1) --(-.7,-.7);
	\draw[color=black,thick] (0.7,0.7) to [out=200, in=90] (.3,0) to [out=270, in=160] (0.7,-0.7);
	\node [below] at (0,-1) {A};
  }}
\newcommand{\makeInB}{\make{
	\draw[color=gray] circle (1);
	\draw[color=black,thick] (0, 1) --(-.7,-.7);
	\draw[color=black,thick] (0,-1) to [out=120, in=180] (.3,-.3) to [out=0, in=90] (0.7,-0.7);
	\draw[color=black,thick] (-.7,.7) to [out=280, in=180] (0,.2) to [out=0, in=260] (0.7,0.7);
	\node [below] at (0,-1) {B};
  }}
\newcommand{\makeInC}{\make{
	\draw[color=gray] circle (1);
	\draw[color=black,thick] (0.7,0.7) to [out=200, in=90] (.3,0) to [out=270, in=160] (0.7,-0.7);
	\draw[color=black,thick] (-0.7,0.7) to [out=275, in=180] (-.3,.3) to [out=0, in=290](0,1);
	\draw[color=black,thick] (-0.7,-0.7) to [out=85, in=180] (-.3,-.3) to [out=0, in=70](0,-1);
	\node [below] at (0,-1) {C};
  }}
\newcommand{\makeInD}{\make{
	\draw[color=gray] circle (1);
	\draw[color=black,thick] (0,1) --(0.7,-0.7);
	\draw[color=black,thick] (-.7,.7) to [out=280, in=180] (0,.2) to [out=0, in=260] (0.7,0.7);
	\draw[color=black,thick] (-0.7,-0.7) to [out=85, in=180] (-.3,-.3) to [out=0, in=70](0,-1);
	\node [below] at (0,-1) {D};
  }}

\title{\textsc{Multi-Skein Invariants for Welded and Extended Welded Knots and Links}}

\author{\textsc{N.\ Backes, M.\ Kaiser, T. Leafblad, E.\ I.\ C.\ Peterson, D.\ N.\ Yetter}}
\date{} 

\begin{document}
\maketitle

\begin{abstract}
    \noindent The theory of welded and extended welded knots is a generalization of classical knot theory.  Welded (resp. extended welded) knot diagrams include virtual crossings (resp. virtual crossings and wen marks) and are equivalent under an extended set of Reidemeister-type moves.
	We present a new class of invariants for welded and extended welded knots and links using a multi-skein relation, following Z. Yang's approach for virtual knots. Using this skein-theoretic approach, we find sufficient conditions on the coefficients to obtain invariance under the extended Reidemeister moves appropriate to welded and extended welded links. 
\end{abstract}

\bigskip

\begin{center}
{\large \textsc{Introduction}
}
\end{center}

\begin{figure}[h]
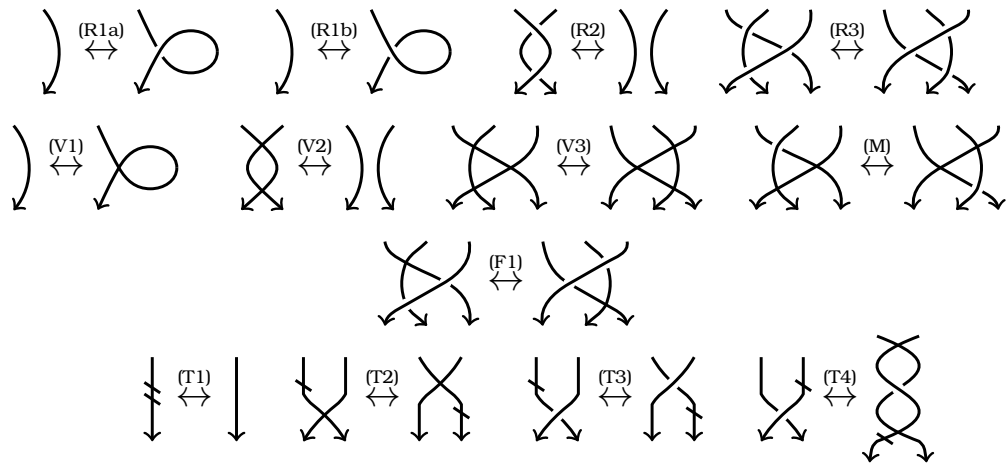

        \centering
    \begin{center}
      \scalebox{1.4}{ \makeRp \makeRn \makeRR \makeRRR}\\
      \scalebox{1.4}{ \makeRv \makeVV \makeVVV \makeM}\\
      \scalebox{1.4}{ \makeW}\\
	  \scalebox{1.4}{ \makeT \makeTt \makeTTT  \makeTTTT}
	\caption{Generalized Reidemeister moves for extended welded links}
	\label{moves} 
    \end{center}
\end{figure}

\noindent The theory of virtual links, introduced by Kauffman in 1999, motivated by the role of Gauss codes in Vassiliev theory and the lack of classical knots corresponding to all possible Gauss codes, consists of the theory of classical links extended by a new type of crossing, the virtual crossing, along with the virtual Reidemeister moves allowing for the description of equivalence classes of virtual link diagrams \cite{K99}. A virtual crossing in a link diagram consists of the simple intersection of two lines, lacking any indication of one strand lying over the other. 
The virtual moves are denoted by V1, V2, and V3, corresponding to the classical Reidemeister moves R1, R2, and R3 respectively. The description  of virtual links is completed by the addition of the mixed move M, in which a classical crossing passes through a pair of virtual crossings in another analogue of R3.  (See Figure \ref{moves}.)

In virtual knot theory the other obvious analogues of R3, in which a virtual crossing passes under or over another strand are absent, and are called ``forbidden moves''.

The theory of Fenn, Rim\'{a}nyi, and Rourke's welded links \cite{FRR} was originally motivated by considering both Artin's braid groups and permutation groups as subgroups of the automorphisms of a free group to define the welded braid group. In this theory, one, and only one, of the forbidden moves, called the welded move is permitted, which allows a virtual crossing to pass under another strand of the link.  Once this move is included, what had been called virtual crossings should be more properly called welded crossings. Following  Damiani \cite{DMarkov}, we denote this move by F1. Several  results on classical links and braids generalize nicely to welded links and braids; in particular, Kauffman and Lambropoulou \cite{KL} proved analogues of the Alexander and Markov theorems for welded links.  Damiani \cite{DLoop} provides an excellent survey of the different geometric descriptions of Fenn, Rim\'{a}nyi, and Rourke's welded braid groups, including as the group of ribbon braids, and as ``loop braids''---motions of an unlink in ${\mathbb R}^3$ in which the components remain in horizontal planes throughout the motion (cf. Baez, Wise and Crans \cite{BWC}).

Finally, one can introduce the wen mark, corresponding in the ribbon or loop pictures to allowing a circle to flip over during the motion, and four moves involving it in order to obtain the theory of extended welded links. The complete family of moves allowed in the theory of extended welded links (and corresponding implicitly to the relations in the extended welded braid groups) are shown in Figure \ref{moves}, where a wen is indicated by an angled mark in the location of the crossing.  Damiani \cite{DMarkov} proved an analogue of Markov's theorem relating the extended welded braid groups to extended welded links in the anticipated way.

The interpretation of welded and extended welded braids as motions of loops in ${\mathbb R}^3$ suggests that surfaces ambient isotopic to the closures of geometric loop braids, with or without wens, will form an interesting class of geometric objects in ${\mathbb R}^4$.  The relevant structure for single components (without wens) -- ribbon torus knots has been known since work of Yajima \cite{Y62} in the early 1960's.

\begin{definition} A {\em ribbon torus knot (resp. link)} is an embedding of the torus $S^1 \times S^1$ (resp. a disjoint union of tori) into ${\mathbb R}^4$, which admits an extension to an immersion of the solid torus $S^1 \times D^2$ (resp. a disjoint union of solid tori) in which the only singularities are disks of double points with a neighborhood modelled by $J:(-1,1) \times D^2 \coprod (-1,1) \times D^2 \rightarrow {\mathbb R}^4$, with $J$ given on the first summand by $(t,x,y) \mapsto (t, 0, x, y)$ and on the second by $(t,x,y) \mapsto (0,t,x/2, y/2)$.
\end{definition}

Satoh \cite{S} proposed that virtual knots could be understood as encoding ribbon torus knots by means of the ``tube map'', but this map is in fact invariant under F1, and thus better seen as an encoding of ribbon torus knots by welded links (cf. the excellent survey by Audoux \cite{A}).  The intial hope that the tube map would be a bijection proved ill-founded, since the horizontal mirror image of a welded link is not in general equivalent under the welded Reidemeister move to the original link, but its image under the tube map is always ambient isotopic in ${\mathbb R}^4$ to the image of the original link.  The introduction of wen marks (cf. Damiani \cite{DMarkov, DLoop}) obviates this problem, as a horizontal mirror image will always be equivalent under extended welded Reidemeister moves to the original link. Moreover, there are other non-welded-equivalent pairs of welded links which become equivalent once it is possible to introduce a pair of wen marks on a component and then move one along the component until they can again be cancelled.   It is, however, unknown at this writing whether this covers all of the sources of non-injectivity for the welded tube map, giving a bijection between extended welded links (with an even number of wens on each component) and ribbon torus links.

\noindent It is the purpose of the present paper to determine the class of invariants of welded and extended welded links determined by multi-skein relations of the form proposed in the context of virtual link theory by Z. Yang \cite{yang}, given in Figure \ref{ZYskein}.

\begin{figure}[h]
    \[ \; \left[ \makeXPD \right]= a \left[\makeXV \right]+ b \left[\makeTT \right]+c \left[\makeC \right] \;\; \]
    \[ \left[ \makeXND \right]= x \left[\makeXV \right]+ y \left[\makeTT \right]+z \left[\makeC \right]  \] 
    \caption{Yang's multi-skein relation}
    \label{ZYskein}

\end{figure}

\begin{center}
    \textsc{\large Construction of Invariants} 
\end{center}

\noindent   Just as Kauffman's skein relation for his bracket polynomial reformulation of the Jones polynomial \cite{KBracket, J} applied to every crossing of an oriented classical knot diagram produces a sum of $2^{\mbox{\small \rm \# crossings}}$ states,  each of which is an unoriented unlink diagram, so applying Yang's relations to the classical crossings of an oriented virtual (resp. welded, extended welded) link diagram produces a sum of $3^{\mbox{\small \rm \# classical crossings}}$ states, each of which is an unoriented virtual (resp. welded, extended welded) link diagram with only virtual (resp. welded, welded) crossings.


It is easy to see that using moves V2 and V3 any state is equivalent to a disjoint union of circles immersed in the plane, each of which can be reduced to an embedded circle by application of some number of instances of V1.
We will not require our Yang-style bracket to be invariant under V1, but will restore V1 invariance by a normalization analogous to the writhe correction used by Kauffman \cite{KBracket}.  We thus introduce three more variables:  a value for removing a component $t$; a value for undoing a virtual crossing to which V1 would apply (a ``virtual kink''), $r$, requiring that $r^2 = 1$ since an analogue of the Whitney trick will allow two virtual kinks on opposite sides of an arc to cancel each other; and a value $s$ for removing a wen, also requiring that $s^2 =1$ to give invariance under T1, giving the skein relations of Figure \ref{TandR}

\begin{figure}[h]
    \begin{center} 
        $\left[ L \cup \makeO\right]=t \left[\;L\;\right]$ \hspace{0.2in}
        $\left[\makeVun \right]= r \left[\makeFun \right] $ \hspace{0.2in}
 	$\left[\makeWun \right]= s \left[\makeFun \right] $\\ 
    \end{center}
    \caption{Skein relations defining $t$, $r$ and $s$}
    \label{TandR}

\end{figure}

As in Kauffman \cite{KBracket} the writhe of a diagram $L$ is given by
\[ w(L) = \mbox{\rm \# positive classical crossings} - \mbox{\rm \# negative classical crossings}.\]  
The virtual writhe $v(L)$ is the number of virtual crossings, and may be taken modulo 2.  Observe that both the writhe and the virtual writhe are preserved by all of the moves in the theory of extended welded links except R1a, R1b, V1 and T4, which (reading from left to right in Figure \ref{moves}) increase the writhe by 1, decrease the writhe by 1, increase the virtual writhe by 1, and decrease the writhe by 2, respectively.

We will proceed by finding conditions on the coefficients $a, b, c, w, y, z, t, r,$ and $s$ which give invariance under the other nine moves, then consider what is required to allow the resulting quantity to be normalized by multiplication by factors depending on the writhe and virtual writhe to obtain invariants of extended welded links.

Of these, invariance under V2, V3, T1, T2 and T3 are immediate without any further constraints on the coefficients.   A moment's consideration shows that satisfying the skein relation and being invariant under V2 and V3 implies invariance under M:  replace the classical crossing by the skein relation, then V3 and V2 (twice) in the summands to which they apply to move the crossing to the other side of the strand passing through both virtual (welded) crossings. 

The general procedure will be to consider the remaining classical and welded Reidemeister moves one by one, imposing relations on coefficients to force our final class of polynomials to be invariants of extended welded links.  In particular, we first consider R2 and F1, as we will see that invariance under these implies invariance under R3. 

By virtue of Damiani's \cite{DMarkov} Markov theorem, it will suffice to verify invariance under the moves in the braid-like position shown in Figure \ref{moves}, there being for instance no need to consider directly the second Reidemeister move with one strand oriented up, the other oriented down.

We begin by imposing invariance under R2. Consider the following expansion of the left-hand side of the diagrammatic representation of R2:\\
\vspace*{-5mm}
        
        

\begin{center}
$\left[ \makeRRhalf \right]= (ay+bx)\left[ \makeXV \right]+ (ax+by) \left[\makeTT\right] +
(azr+cxr+bz+cy+czt)\left[\makeC\right] $.
\end{center}

\noindent To have invariance under this move, we must have
\vspace*{-1mm}
\begin{align*}
    &ay+bx=0, \\
    &ax+by=1, \\
    &azr+cxr+bz+cy+czt=0,
\end{align*} 

\noindent which yields the relations
\vspace*{-1mm}
\begin{align*}
x&=\frac{-a}{b^2-a^2},\\
y&=\frac{b}{b^2-a^2},\\
z&=\frac{rac-bc}{(b^2-a^2)(ra+b+ct)}.
\end{align*}\\
\vspace*{-5mm}

Expanding F1 in a similar manner does not produce terms that allow us to use this method. Instead, having expanded both sides of F1 and reduced terms by pulling out factors of $t$ and applying any of V1, V2, V3, and R2 as needed, and cancelling like terms, we obtain the following equation:

\begin{center}
    \includegraphics[height=8.5mm,scale=0.5]{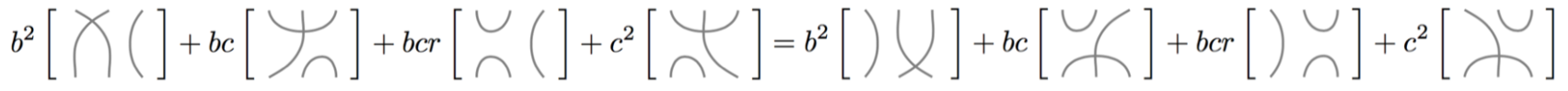}
\end{center}

It is clear that we cannot guarantee equality with any nontrivial choice of coefficients. However, it is sufficient to show that equality (in particular, equality of the four mismatched terms) is restored for each way of closing the depicted tangles to form local states (unoriented links with only welded crossings). This can be verified by checking the fifteen ways in which the six endpoints are connected by strands with only virtual crossings. 


\begin{figure}[h]
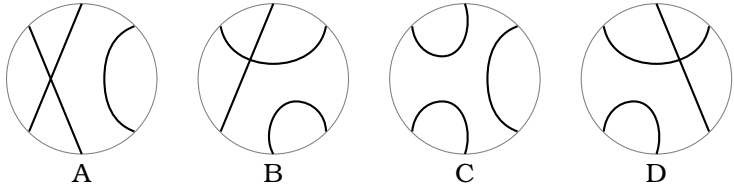

        \centering 
    	\makeInA \makeInB \makeInC \makeInD \\
     \caption{The mismatched tangles in the expansion of F1. To obtain the other set, rotate each by 180 degrees.}
     \label{insidez} 
\end{figure}

\begin{figure}[h]
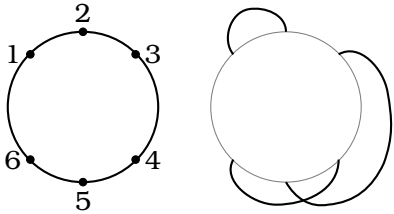

        \centering 
    	\makeCases \makecB
     \caption{Example of a closure in which 1 and 2, 3 and 5, and 4 and 6 are connected.}
     \label{cases_exe}
\end{figure}

To systematically verify these cases, we label each endpoint as in Figure \ref{cases_exe} and refer to a case by writing three pairs of endpoint labels indicating which endpoints are connected. For example, we will denote the case given by the closure shown in Figure \ref{cases_exe} by 12:35:46.

Computing the equations arising in each of the 15 cases yields 3 nontrivial equations

\vspace*{-4mm}
\begin{align*}
	& b^2+bc+bct+c^2=b^2t+bct^2+bc+c^2t, \\
	& b^2+2bct+c^2t^2=b^2t+2bc+c^2, and\\
	& b^2t^2+2bct+c^2=b^2+2bc+c^2t,
\end{align*}

\noindent corresponding to the closures 12:35:46 (13:26:45), 13:24:56 (15:23:46), and 15:26:34 (16:24:35), respectively

It follows that invariance under F1 requires either 
\vspace*{-1mm}
\begin{align*}
b&=\pm c  \hspace{.1in} t = \mp2 \text{ or}\\
t&=1.
\end{align*}

In the case $t=1$, it is easy to verify that any link evaluates to $1$.  In what follows we will consider only the case $b=\pm c$ with $t = \mp2$.

One might hope that applying the same procedure to the before and after states of R2 might yield more solutions, but this is not the case.

Invariance of the bracket under R2 and F1 gives invariance under R3 by an argument analogous to that by which V2 and V3 gave invariance under M:  apply the skein relation to the crossing at the bottom left of the left-hand side of R3, then use R2 and F1 to move its remants in the summands over the other strand, and reverse the skein relation to give the right-hand side of R3.

We have now accomplished our first goal of finding conditions on the coefficients that give invariance under all of the moves except R1a, R1b, V1 and T4.  It turns out that further calculations can be simplified by letting $
\nu = c/b$ and then letting $c = \nu b$ and $t = -2\nu$, with $\nu^2 =1$.

Substituting $\nu b$ for $c$ and $-2\nu$ for $t$ in the conditions expressing $x, y$ and $z$ in terms of $a, b$ and $c$ to give R2 invariance simplies them to

\begin{align*}
        x&=\frac{-a}{b^2-a^2},\\
        y&=\frac{b}{b^2-a^2},\\
        z&=\frac{\nu b}{b^2-a^2}.
        \end{align*}

\bigskip

For brevity in what follows, we let $\delta = b^2 - a^2$, letting us rewrite the preceding as
$x = -a/\delta$, $y = b/\delta$ and $z = \nu b/\delta$.

It is easy to see that invariance under V1 can be restored by multiplying the bracket by $r^{v(L)}$.  To restore invariance under R1, observe that applying the skein relation to the crossing in the right hand side of R1a and R1b and simplifying the result using the relations involving $t$ and $r$ and those just given gives

\begin{center}
    $\left[\makePun\right] = (ar - \nu b)\left[\makeFunn \right]
    \;\;\text{and}\;\; \left[\makeNun\right] = (\frac{-ra - \nu b}{b^2 - a^2})\left[\makeFunn \right].$ 
\end{center}

\noindent Observe that the two coefficients in the right-hand sides of these equations are reciprocals (assuming 
$a^2 \neq b^2$ or equivalently $\delta \neq 0$), so we may restore invariance under R1a and R1b by multiplying by $(ar - \nu b)^{-w(L)}$.

We thus have

\begin{theorem}
The quantity

\begin{center}
    $Y(L) = r^{v(L)}(ar - \nu b)^{-w(L)}[L]$
\end{center} 

\noindent is an invariant of welded links, where $[L]$ denotes the rational function of $a, b$ and the $\pm 1$-valued variables $r$ and $\nu$ given by the skein relations of Figure \ref{goodwelded}.

\end{theorem}

\begin{figure}[h]
    \[ \; \left[ \makeXPD \right]= a \left[\makeXV \right]+ b \left[\makeTT \right]+ \nu b \left[\makeC \right] \;\; \]
    \[ \left[ \makeXND \right]= \frac{-a}{\delta} \left[\makeXV \right]+ \frac{b}{\delta}  \left[\makeTT \right] + \frac{\nu b}{\delta}  \left[\makeC \right]  \] 
\[ \left[ L \cup \makeO\right]= -2\nu \left[\;L\;\right] \] 
 \[ \left[\makeVun \right]= r \left[\makeFun \right] \] 

    \caption{Skein relations for welded links}
    \label{goodwelded}

\end{figure}

In fact, provided $\nu =1$, the same conditions give an invariant of extended welded links.

For brevity, let $\omega = ar - \nu b$ and thus $\omega^{-1} = -ar/\delta - \nu b /\delta$.

For T4 to hold for the writhe-corrected bracket of the previous theorem, we must have the three closures of

\[ a \left[\makeXV \right]+ b \left[\makeTT \right]+ \nu b \left[\makeC \right] \]

\noindent evaluate to exactly $\omega^2$ times the corresponding closure of

\[ -a /\delta \left[\makeXV \right]+ b/\delta \left[\makeTT \right]+ \nu b/\delta \left[\makeC \right] \]

In each case, multiplying the equation by $\delta$ results in 
\begin{align*}
        \delta(ar - b ) &=\omega^2(-ra - b),\\
        \delta(ar - \nu b ) &=\omega^2(-ra - \nu b),\\
        \delta(-2 a \nu + r b + \nu r b ) &=\omega^2(2\nu a + r b + \nu r b).
        \end{align*}

If $\nu = -1$, the last equation implies that either $a$ or $\delta + \omega^2$ must be zero.

In the latter case, either $b$ is zero, and the invariant virtualizes all crossings, and thus, depends only on the number of components and the number of wens, and will simply be 
\[ s^{\mbox{\rm \# of wens}} 2^{\mbox{\rm \# of components of $L$}}, \] 
\noindent or $\omega^2 - \delta$ is zero, either contradicting the condition $\delta \neq 0$ or reducing to the case $\nu = 1$ in charactistic 2 with $\delta + \omega^2 = 0$ (which holds vaccuously in characteristic 2).

If $a = 0$,  the invariant again reduces to $s^{\mbox{\rm \# of wens}} 2^{\mbox{\rm \# of components of $L$}}$, although the proof of this is less direct than that for $b =0$.

However, if $\nu = 1$, all three of the equations simplify to

\[ (\delta + \omega^2) r a + (\omega^2 - \delta) b = 0. \]

\noindent Recalling that $\delta = b^2 - a^2$ and $\omega = ar - \nu b$, here $\omega = ra - b$, and substituting, we find that this equation holds for any values of $a$ and $b$, thus showing

\begin{theorem}
The quantity

\begin{center}
    $Y(L) = r^{v(L)}(ar -  b)^{-w(L)}[L]$
\end{center} 

\noindent is an invariant of extended welded links where $[L]$ denotes the rational function of $a, b$ and the $\pm 1$-valued variables $r$ and $s$ given by the skein relations of Figure \ref{goodextwelded}.

\end{theorem}

\begin{figure}[h]
    \[ \; \left[ \makeXPD \right]= a \left[\makeXV \right]+ b \left[\makeTT \right]+ b \left[\makeC \right] \;\; \]
    \[ \left[ \makeXND \right]= \frac{-a}{\delta} \left[\makeXV \right]+ \frac{b}{\delta}  \left[\makeTT \right] + \frac{b}{\delta}  \left[\makeC \right]  \] 
\[ \left[ L \cup \makeO\right]= -2\left[\;L\;\right] \] 
 \[ \left[\makeVun \right]= r \left[\makeFun \right] \] 
\[ \left[\makeWun \right]= s \left[\makeFun \right] \]
    \caption{Skein relations for extended welded links ($\nu = 1$)}
    \label{goodextwelded}

\end{figure}

We conclude by showing that this invariant of extended welded links is non-trivial.

Plainly 

\[ Y(\mbox{\rm 2 component unlink}) = 4. \] 

\noindent However, simple calculations shows that 

\[ Y(\mbox{\rm positive Hopf link}) = \frac{(4a^2 + 4rab + 12b^2)(a^2 + 2rab + b^2)}{(b^2 - a^2)^2}. \]

\noindent While if $L$ is the ``one-virtual, one-positive Hopf link'' 

\[ Y(L) = \frac{4a^2 - 4rab}{b^2 - a^2}. \]

\begin{center}
    \textsc{\large Prospects for Future Research} 
\end{center}

One question which immediately arises for any invariant of virtual, welded, or extended welded links is whether its restriction to classical links gives anything new.  In the sense of the naive question ``Does $Y(L)$ coincide with any of the classical knot polynomials (with Zariski dense sets of values) coming from quantum groups?" the negative answer suggests it is something new.  However, we conjecture that a more strict interpretion of the question phrased oppositely: ``Does $Y(L)$ distinguish any pairs of classical links not distinguished by any of the invariants coming from quantum groups?" also has a negative answer.  We would propose studying this question by substituting formal power series for the coefficients in $Y(L)$ and examining the resulting sequences of Vassiliev invariants.

It should be pointed out that before this is done, that it might be worthwhile doing a change of variables.  $Y(L)$ as described thus far is a rational function of $a$ and $b$ (also depending on the $\pm 1$ valued variables $r, \nu$ and $s$), but the denominators are always powers of $\delta = b^2 - a^2$.  If the change of variables $\alpha := b+a$, $\beta := b-a$ is applied in characteristic different from 2,  meaning $b = \frac{1}{2}(\alpha + \beta)$ and $a = \frac{1}{2}(\alpha - \beta)$,  the rational function $Y(L)$ reduces to a Laurent polynomial in $\alpha$ and $\beta$, homogeneous of degree 0, which can then be dehomogenized to give a one-variable Laurent polynomial.
It would probably be most fruitful to substitute a formal power series (say $\exp(x)$) into this one-variable polynomial and to compare the power of the resulting sequence of Vassiliev invariants with those arising from invariants coming from quantum groups.

Once the question of whether the inclusion of wens suffices to explain all of the non-injectivity of the tube map has been resolved, should there be additional equivalences required, it will be important to examine the behavior of $Y(L)$ within the equivalence classes (ideally under the additional moves, if the needed equivalence is induced by move).  Should the tube map prove to be injective on extended welded links, of course, this would imply that $Y(L)$ in the case $\nu = 1$ is an invariant of ribbon torus links.

It would also be worthwhile to investigate other skein relations in the context of extended welded links.

\newpage
\noindent \textsc{\Large Acknowledgements}
\newline

\noindent The authors thank Kansas State University for hosting the SUMaR REU. This research was supported by National Science Foundation Award \# DMS-1659123.
\newline

\textsc{University of Minnesota, Minneapolis, MN} (Current Address of N.\ Backes)

\textsc{North Central College, Naperville, IL} (Current Address of M.\ Kaiser)

\textsc{Truman State University, Kirksville, MO} (Current Address of T.\ Leafblad)

\textsc{University of Minnesota, Minneapolis, MN} (Current Address of E.\ I.\ C.\ Peterson)

\textsc{Kansas State University, Manhattan, KS} (Current Address of D.\ N.\ Yetter)

\newpage
\renewcommand\refname{\normalfont{\textsc{References}}}

\end{document}